%% file: BeukersStyle.tex
\input a4

\input library.tex

\allowdisplaybreaks[4]

\begin {document}

\begin{center}
{\bf \scshape Generalizations of Guillera-Sondow's double integral formulas} \\
{\bf \scshape Sergey Zlobin}
\end{center}

{\small
{\it Abstract.}
We evaluate certain multidimensional integrals in terms
of the Lerch transcendent function $\Phi$, generalizing Guillera-Sondow's
formulas. As an application, we get new representations of classical constants
like Euler's constant $\gamma$ and $\ln(4/\pi)$.
}

\vspace{5mm}

The {\it Lerch transcendent} $\Phi$ is defined as
the analytic continuation of the series
\begin{equation*}
\Phi(z,s,u)=\frac{1}{u^s}+\frac{z}{(u+1)^s}+\frac{z^2}{(u+2)^s}+\cdots,
\end{equation*}
which converges for any complex number $u$ with $\Re u>0$
if $z$ and $s$ are any
complex numbers with either $|z|<1$, or $|z|=1$ and $\Re s > 1$
(we suppose $\zeta^s=\exp(s \log \zeta)$,
where $\log \zeta$ is the principal branch of the logarithm).
The function $\Phi$ is holomorphic in $z$ and $s$, 
for $z \in \C \backslash [1,\infty]$
and all complex $s$ (see \cite[section 1.11]{Bateman} or 
\cite[section 2]{GuilleraSondow}).
From the definition it follows that
\begin{equation}
\label{LerchEq}
\Phi(z,s,u+1) = \frac{1}{z} \left(  \Phi(z,s,u)- \frac{1}{u^s} \right),
\end{equation}
\begin{equation}
\label{LerchEq2}
\Phi(z,s+1,u) = -\frac{1}{s} \frac{\partial \Phi}{\partial u}(z,s,u).
\end{equation}

In the paper \cite{GuilleraSondow} J. Sondow and J. Guillera proved 
the following two theorems.
\begin{theorem}
\label{SondowGuilleraTh1}
Suppose $u>0$, $v>0$, $u \ne v$,
either $z \in \C \backslash [1,\infty)$ and $\Re s > -2$, or $z=1$ and $\Re s > -1$.
Then
\begin{align*}
\int_{\cube{2}} \frac{x_1^{u-1} x_2^{v-1}}{1-z x_1 x_2}
(-\ln x_1 x_2)^{s} dx_1 dx_2 &= \Gamma(s+1)
\frac{\Phi(z,s+1,v)-\Phi(z,s+1,u)}{u-v}, \\
\int_{\cube{2}} \frac{(x_1 x_2)^{u-1}}{1-z x_1 x_2}
(-\ln x_1 x_2)^{s} dx_1 dx_2 &= {\Gamma(s+2)} \Phi(z,s+2,u).
\end{align*}
\end{theorem}
\begin{theorem}
\label{SondowGuilleraTh2}
Suppose $u>0$, 
either $z \in \C \backslash  [1, \infty)$ and $\Re s > -3$, or $z=1$ and $\Re s > -2$.
Then
\begin{multline*} 
\int_{\cube{2}} \frac{1-x_1}{1-z x_1 x_2} (x_1 x_2)^{u-1}
(-\ln x_1 x_2)^{s} dx_1 dx_2
\\ =\Gamma(s+2) \left[ \Phi(z,s+2,u)+\frac{(1-z)
\Phi(z,s+1,u)-u^{-s-1}}{z (s+1)} \right].
\end{multline*}
\end{theorem}
The purpose of this paper is to prove the following $m$-dimensional analogs
of Theorems \ref{SondowGuilleraTh1} and \ref{SondowGuilleraTh2}
(in what follows $d \overline x$ means $dx_1 dx_2 \cdots dx_m$, where
$m$ is the dimension of an integral).
\begin{theorem}
\label{MainTh1}
Suppose $m$ is a positive integer, $\Re u>0$, $\Re v>0$, $u \ne v$,
either $z \in \C \backslash [1,\infty)$ and $\Re s > -m$,
or $z=1$ and $ \Re s > 1-m$. For the case $m>1$ define the function
\begin{multline*}
F_{m,u,v}(x_1,x_2,\dots,x_m) =
(x_1 x_2 \cdots x_m)^{v-1} (x_1^{u-v} + (x_1 x_2)^{u-v} + \cdots + 
(x_1 x_2 \cdots x_{m-1})^{u-v} ).
\end{multline*}
Then
\begin{multline}
\int_{\cube{m}} \frac{F_{m,u,v}(x_1,x_2,\dots,x_m)}
{1-z x_1 x_2 \cdots x_m}
(-\ln x_1 x_2 \cdots x_m)^{s} d \overline x \\
=\frac{\Gamma(s+m-1)}{(m-2)!} \cdot
\frac{\Phi(z,s+m-1,v)-\Phi(z,s+m-1,u)}{u-v} \quad \mbox{for } m>1
\label{MainTh1Eq1},
\end{multline}
\begin{equation}
\label{MainTh1Eq2}
\int_{\cube{m}} \frac{(x_1 x_2 \cdots x_m)^{u-1}}{1-z x_1 x_2 \cdots x_m}
(-\ln x_1 x_2 \cdots s_m)^{s} d \overline x =
\frac{\Gamma(s+m)}{(m-1)!} \Phi(z,s+m,u) \quad \mbox{for } m \ge 1.
\end{equation}
\end{theorem}
\begin{theorem}
\label{MainTh2}
Suppose $m$ is an integer $>1$, $\Re u>0$,
either $z \in \C \backslash [1, \infty)$ and $\Re s > -m-1$, or $z=1$ and $\Re s > -m$.
Then
\begin{multline} 
\int_{\cube{m}} \frac{m-1 - x_1-x_1 x_2 - \cdots - x_1 x_2 \cdots x_{m-1}}
{1-z x_1 x_2 \cdots x_m} (x_1 x_2 \cdots x_m)^{u-1}
(-\ln x_1 x_2 \cdots x_m)^{s} \, d \overline x
\\ = \frac{\Gamma(s+m)}{(m-2)!} \left[ \Phi(z,s+m,u)+\frac{(1-z)
\Phi(z,s+m-1,u)-u^{-s-m+1}}{z (s+m-1)} \right]
\label{MainTh2Eq}.
\end{multline}
\end{theorem}

In the case $m=2$ Theorems \ref{MainTh1} and \ref{MainTh2} give
Theorems \ref{SondowGuilleraTh1} and \ref{SondowGuilleraTh2}.
As an example, we give also the case $m=3$.
\begin{example} 
a)
If $\Re u>0$, $\Re v>0$, $u \ne v$,
either $z \in \C \backslash [1,\infty)$ and $\Re s > -3$,
or $z=1$ and $ \Re s > -2$, then
$$
\int_{\cube{3}}
\frac{x_1^{u-1} x_2^{v-1} x_3^{v-1} + x_1^{u-1} x_2^{u-1} x_3^{v-1}}
{1-z x_1 x_2 x_3}
(-\ln x_1 x_2 x_3)^{s} d \overline x \\
=\Gamma(s+2)
\frac{\Phi(z,s+2,v)-\Phi(z,s+2,u)}{u-v}.
$$
b) If $\Re u>0$, 
either $z \in \C \backslash [1,\infty)$ and $\Re s > -3$,
or $z=1$ and $ \Re s > -2$, then
$$
\int_{\cube{3}}
\frac{(x_1 x_2 x_3)^{u-1}}
{1-z x_1 x_2 x_3}
(-\ln x_1 x_2 x_3)^{s} d \overline x \\
=\frac{\Gamma(s+3)}{2}
\Phi(z,s+3,u).
$$
c) If $\Re u>0$, 
either $z \in \C \backslash [1,\infty)$ and $\Re s > -4$,
or $z=1$ and $ \Re s > -3$, then
\begin{multline*}
\int_{\cube{3}} \frac{2 - x_1-x_1 x_2}
{1-z x_1 x_2 x_3} (x_1 x_2 x_3)^{u-1}
(-\ln x_1 x_2 x_3)^{s} \, d \overline x
\\ = \Gamma(s+3) \left[ \Phi(z,s+3,u)+\frac{(1-z)
\Phi(z,s+2,u)-u^{-s-2}}{z (s+2)} \right].
\end{multline*}
\end{example}

In \cite{GuilleraSondow} many interesting applications
of Theorems \ref{SondowGuilleraTh1} and \ref{SondowGuilleraTh2} are given.
All of them can be generalized by Theorems
\ref{MainTh1} and \ref{MainTh2}; indeed, by these four theorems we have
\begin{multline*}
\int_{\cube{2}} \frac{x_1^{u-1} x_2^{v-1}}{1-z x_1 x_2}
(-\ln x_1 x_2)^{s} dx_1 dx_2 \\
= (m-2)!
\int_{\cube{m}} \frac{F_{m,u,v}(x_1,x_2,\dots,x_m)}
{1-z x_1 x_2 \cdots x_m}
(-\ln x_1 x_2 \cdots x_m)^{s-m+2} d \overline x
\quad \mbox{for } m>1,
\end{multline*}
\begin{multline*}
\int_{\cube{2}} \frac{(x_1 x_2)^{u-1}}{1-z x_1 x_2}
(-\ln x_1 x_2)^{s} dx_1 dx_2 \\
= (m-1)!
\int_{\cube{m}} \frac{(x_1 x_2 \cdots x_m)^{u-1}}
{1-z x_1 x_2 \cdots x_m}
(-\ln x_1 x_2 \cdots x_m)^{s-m+2} d \overline x
\quad \mbox{for } m \ge 1,
\end{multline*}
\begin{multline*}
\int_{\cube{2}} \frac{1-x_1}{1-z x_1 x_2} (x_1 x_2)^{u-1}
(-\ln x_1 x_2)^{s} dx_1 dx_2 \\
=(m-2)!
\int_{\cube{m}} \frac{m-1 - x_1-x_1 x_2 - \cdots - x_1 x_2 \cdots x_{m-1}}
{1-z x_1 x_2 \cdots x_m} (x_1 x_2 \cdots x_m)^{u-1} \\
\times (-\ln x_1 x_2 \cdots x_m)^{s-m+2} \, d \overline x
\quad \mbox{for } m>1.
\end{multline*}

We give here only two examples.
\begin{example}
Let m be an integer $>1$, and $\gamma=\lim_{n \to \infty}
\left( 1+\frac{1}{2}+\cdots+\frac{1}{n} - \ln n \right)$ be
Euler's constant. Then
$$
\gamma=(m-2)! \int_{\cube{m}}
\frac{m-1 - x_1-x_1 x_2 - \cdots - x_1 x_2 \cdots x_{m-1}}
{(1-x_1 x_2 \cdots x_m) (-\ln x_1 x_2 \cdots x_m)^{m-1}}
d \overline x.
$$
\end{example}
\begin{example}
For an integer $m>1$ the following identity holds:
$$
\ln \frac{4}{\pi} = (m-2)! \int_{\cube{m}}
\frac{m-1 - x_1-x_1 x_2 - \cdots - x_1 x_2 \cdots x_{m-1}}
{(1+x_1 x_2 \cdots x_m) (-\ln x_1 x_2 \cdots x_m)^{m-1}}
d \overline x.
$$
\end{example}
We omit details of proofs of these examples and refer to the case 
$m=2$, which was considered by J. Sondow \cite{SondowMonthly}.

To prove Theorem \ref{MainTh1}, we will require two lemmas.
The first is the identity (\ref{MainTh1Eq2}) for $m=1$,
and is classical.
\begin{lemma}
\label{Lemma1}
Suppose $\Re u>0$,
either $z \in \C \backslash [1,\infty)$ and $\Re s > -1$,
or $z=1$ and $ \Re s > 0$. Then
$$
\int_0^1 \frac{x^{u-1}}{1-z x} (-\ln x)^{s} dx =
\Gamma(s+1) \Phi(z,s+1,u).
$$
\end{lemma}
\proof
The integral, call it $I$, defines a holomorphic function of $z$ and $s$
under the conditions
stated. We can prove the statement for $|z|<1$ and $\Re s > 0$ and then
use analytic continuation. Expand $1/(1-zx)$ into a geometric series and
then integrate:
$$
I=\sum_{n=0}^{\infty} z^n \int_0^1 x^{u+n-1} (-\ln x)^{s} dx.
$$
Making the substitution $x=e^{-y}$, we obtain
$$
I
=\sum_{n=0}^{\infty} z^n \int_0^{\infty} e^{-(u+n)y} y^{s} dy
=\sum_{n=0}^{\infty} \frac{\Gamma(s+1) z^n}{(u+n)^{s+1}}
=\Gamma(s+1) \Phi(z,s+1,u),
$$
and the lemma follows.

\begin{lemma}
Let $\alpha \ne 0$ and $x \in (0,1]$.
Then the following identities hold for $k \ge 1$: \\
a) 
\begin{equation}
\label{SimplexEq1}
\int_{1 \ge t_1 \ge t_2 \ge \cdots \ge t_k \ge x}
\frac{1}{t_1 t_2 \cdots t_k} dt_1 dt_2 \cdots dt_k = \frac{(-\ln x)^k}{k!},
\end{equation}
b)
\begin{equation}
\label{SimplexEq2}
\int_{1 \ge t_1 \ge t_2 \ge \cdots \ge t_k \ge x}
\frac{t_1^{\alpha} + t_2^{\alpha} + \cdots + t_k^{\alpha}}
{t_1 t_2 \cdots t_k} dt_1 dt_2 \cdots dt_k
= \frac{(-\ln x)^{k-1}}{(k-1)!} \cdot \frac{1-x^\alpha}{\alpha}. 
\end{equation}
\end{lemma}
\proof 
The identity (\ref{SimplexEq1}) is easily 
proved using induction and the equality
$$
\int_{1 \ge t_1 \ge t_2 \ge \cdots \ge t_k \ge x}
\frac{dt_1 dt_2 \cdots dt_k}{t_1 t_2 \cdots t_k} =
\int_x^1 \left( 
\int_{1 \ge t_1 \ge t_2 \ge \cdots \ge t_{k-1} \ge t_k}
\frac{dt_1 dt_2 \cdots dt_{k-1}}{t_1 t_2 \cdots t_{k-1}} 
\right) \frac{dt_k}{t_k}.
$$
Denote the integral in (\ref{SimplexEq2}) by $I_{k}(x)$. We prove by
induction; the case $k=1$ is true. Suppose $k>1$ and
the statement is true for $k-1$.
We have
$$
I_{k}(x)=
\int_x^1 I_{k-1}(t_{k}) \frac{dt_k}{t_k} + 
\int_x^1 \left( 
\int_{1 \ge t_1 \ge t_2 \ge \cdots \ge t_{k-1} \ge t_k}
\frac{dt_1 dt_2 \cdots dt_{k-1}}{t_1 t_2 \cdots t_{k-1}} 
\right) t_k^{\alpha-1} dt_k.
$$
Apply (\ref{SimplexEq1}) to the integral in parentheses:
$$
I_{k}(x)=
\int_x^1 I_{k-1}(t_{k}) \frac{dt_k}{t_k} + 
\int_x^1 \frac{(-\ln t_k)^{k-1}}{(k-1)!} t_k^{\alpha-1} dt_k.
$$
Using the induction hypothesis, we obtain
\begin{align*}
I_{k}(x)=&
\int_x^1
\left( 
\frac{(-\ln t_k)^{k-2}}{(k-2)!} \cdot \frac{1-t_k^\alpha}{\alpha t_k}
+ \frac{(-\ln t_k)^{k-1}}{(k-1)!} t_k^{\alpha-1}
\right) dt_k \\
=& \left. \frac{(-\ln t_k)^{k-1}}{(k-1)!} \cdot \frac{t_k^{\alpha}-1}{\alpha}
\right|_x^1
=\frac{(-\ln x)^{k-1}}{(k-1)!} \cdot \frac{1-x^{\alpha}}{\alpha}.
\end{align*}
Now the lemma is completely proved.

\vspace{2mm} \noindent
{\bf Proof of Theorem \ref{MainTh1}.}
The integrals $J_1$ and $J_2$ in (\ref{MainTh1Eq1}) and (\ref{MainTh1Eq2})
define holomorphic functions of $s$ under the conditions
stated. We can prove the statement for $\Re s > 0$ and then
use analytic continuation.

First we prove (\ref{MainTh1Eq2}). Make the substitution
\begin{equation}
\label{Subst}
x_1=t_1, \quad x_2=t_2/t_1, \quad x_3=t_3/t_2, 
\quad \dots, \quad x_m=t_m/t_{m-1}
\end{equation}
in $J_2$. We obtain
\begin{align*}
J_2= & \int_{1 \ge t_1 \ge t_2 \ge \cdots \ge t_m \ge 0}
\frac{t_m^{u-1}}{1-z t_m} (-\ln t_m)^s \frac{1}{t_1 t_2 \cdots t_{m-1}}
dt_1 dt_2 \cdots dt_m \\
= & \int_{0}^1 \frac{t_m^{u-1}}{1-z t_m} (-\ln t_m)^s
\left(
\int_{1 \ge t_1 \ge t_2 \ge \cdots \ge t_{m-1} \ge t_m}
\frac{1}{t_1 t_2 \cdots t_{m-1}}
dt_1 dt_2 \cdots dt_{m-1} \right)
dt_m.
\end{align*}
Applying (\ref{SimplexEq1}) with $x=t_m$ and $k=m-1$, we get
$$
J_2= \frac{1}{(m-1)!}
\int_{0}^1 \frac{t_m^{u-1}}{1-z t_m} (-\ln t_m)^{s+m-1} dt_m.
$$
It remains to apply Lemma \ref{Lemma1}.

Now we prove (\ref{MainTh1Eq1}). Denote $\alpha=u-v$, then
$$
J_1=\int_{\cube{m}} \frac{(x_1 x_2 \cdots x_m)^{v-1}}{1-z x_1 x_2 \cdots x_m}
(x_1^{\alpha} + (x_1 x_2)^{\alpha} + \cdots + (x_1 x_2 \cdots x_{m-1})^{\alpha})
(-\ln x_1 x_2 \cdots s_m)^{s} d \overline x.
$$
Make the substitution (\ref{Subst})
$$
J_1=
\int_{0}^1 \frac{t_m^{v-1}}{1-z t_m} (-\ln t_m)^s
\left(
\int_{1 \ge t_1 \ge t_2 \ge \cdots \ge t_{m-1} \ge t_m}
\frac{t_1^{\alpha} + t_2^{\alpha} + \cdots + t_{m-1}^{\alpha}}
{t_1 t_2 \cdots t_{m-1}}
dt_1 dt_2 \cdots dt_{m-1} \right)
dt_m
$$
and apply (\ref{SimplexEq1})
$$
J_1=\frac{1}{(m-2)! \alpha} \left(
\int_{0}^1 \frac{t_m^{v-1}}{1-z t_m} (-\ln t_m)^{s+m-2} dt_m -
\int_{0}^1 \frac{t_m^{v+\alpha-1}}{1-z t_m} (-\ln t_m)^{s+m-2} dt_m
\right).
$$
It remains to apply Lemma \ref{Lemma1} to both integrals 
and get back to $u$ from $\alpha$. The theorem is proved. \\
{\bf Remark.} The formula (\ref{MainTh1Eq2}) can be also obtained
by letting $v \to u$ in (\ref{MainTh1Eq1}) and using the identity
(\ref{LerchEq2}).

\vspace{2mm} \noindent
{\bf Proof of Theorem \ref{MainTh2}.}
The integral $J$ in (\ref{MainTh2Eq})
defines a function which is holomorphic in $s$,
when $\Re s > -m-1$ if $z \in \C \backslash [1,\infty]$, and when
$\Re s > -m $ if $z=1$.
We prove the statement for $\Re s > 0$ and then use analytic continuation.
We have
\begin{multline*}
J=(m-1) 
\int_{\cube{m}} \frac{(x_1 x_2 \cdots x_m)^{u-1}}
{1-z x_1 x_2 \cdots x_m}
(-\ln x_1 x_2 \cdots x_m)^{s} d \overline x \\
-
\int_{\cube{m}} \frac{F_{m,u+1,u}(x_1, x_2, \dots, x_m)}
{1-z x_1 x_2 \cdots x_m} (-\ln x_1 x_2 \cdots x_m)^{s} \, d \overline x.
\end{multline*}
Apply Theorem \ref{MainTh1} to both integrals:
\begin{align*}
J= & (m-1) \frac{\Gamma(s+m)}{(m-1)!} \Phi(z,s+m,u) \\
& -
\frac{\Gamma(s+m-1)}{(m-2)!} (\Phi(z,s+m-1,u)-\Phi(z,s+m-1,u+1)) \\
= & \frac{\Gamma(s+m)}{(m-2)!} \left[ 
\Phi(z,s+m,u) + \frac{\Phi(z,s+m-1,u+1)-\Phi(z,s+m-1,u)}{(s+m-1)}
\right].
\end{align*}
Use (\ref{LerchEq}) and the theorem follows.

\vspace{3mm}

The way which we prove Theorem \ref{MainTh1} can be applied to any
integral
$$
\int_{\cube{m}} \frac{x_1^{u_1} x_2^{u_2} \cdots x_m^{u_m}}
{1-z x_1 x_2 \cdots x_m}
(-\ln x_1 x_2 \cdots x_m)^{s} d \overline x.
$$
We give the formula for the case when all $u_i$ are different.
\begin{theorem}
\label{MainTh3}
Suppose $m \ge 1$ and $\Re u_1>0$ , $\Re u_2>0$, \dots, $\Re u_m>0$,
$u_i \ne u_j$ whenever $i \ne j$, and
either $z \in \C \backslash [1,\infty)$ and $\Re s > -1$,
or $z=1$ and $ \Re s > 0$. Then the following identity holds:
\begin{equation}
\label{MainTh3Eq}
\int_{\cube{m}} \frac{x_1^{u_1-1} x_2^{u_2-1} \cdots x_m^{u_m-1}}
{1-z x_1 x_2 \cdots x_m}
(-\ln x_1 x_2 \cdots x_m)^{s} d \overline x=
\Gamma(s+1) \sum_{i=1}^m \frac{\Phi(z,s+1,u_i)}{\prod_{j \ne i} (u_j-u_i)}.
\end{equation}
\end{theorem}
To prove Theorem \ref{MainTh3} we require the following
\begin{lemma}
Let $k \ge 1$ and $u_1$, $u_2$, \dots, $u_{k+1}$ be arbitrary numbers with
$u_i \ne u_j$ whenever $i \ne j$, and $x \in (0,1]$.
Then the following identity hold: \\
\begin{equation}
\label{SimplexEq3}
\int_{1 \ge t_1 \ge t_2 \ge \cdots \ge t_k \ge x}
t_1^{u_1-u_2-1} t_2^{u_2-u_3-1} \cdots t_k^{u_k-u_{k+1}-1}
dt_1 dt_2 \cdots dt_k =
\sum_{i=1}^{k+1} \frac{x^{u_i-u_{k+1}}}{\prod_{j=1, j \ne i}^{k+1} (u_j-u_i)}.
\end{equation}
\end{lemma}
\proof
Denote the integral in (\ref{SimplexEq3}) by $I(u_1,u_2,\dots,u_{k+1};x)$.
We prove by induction; the case $k=1$ is true. Suppose $k>1$ and
the statement is true for $k-1$, then
\begin{align*}
I(u_1,u_2,\dots,u_{k+1};x)=&
\int_x^1 I(u_1,u_2,\dots,u_{k};t_k) t_k^{u_k-u_{k+1}-1} dt_k \\
=&\int_x^1 \sum_{i=1}^{k} \frac{t_k^{u_i-u_{k}}}{\prod_{j=1, j \ne i}^k (u_j-u_i)}
t_k^{u_k-u_{k+1}-1} dt_k \\
=& \sum_{i=1}^{k} \frac{1}{\prod_{j=1, j \ne i}^k (u_j-u_i)} 
\int_x^1 t_k^{u_i-u_{k+1}-1} dt_k \\
= & \sum_{i=1}^{k} \frac{1}{\prod_{j=1, j \ne i}^k (u_j-u_i)} \cdot
\frac{ 1 - x^{u_{i}-u_{k+1}} }{u_i-u_{k+1}} \\
= & \sum_{i=1}^{k} \frac{x^{u_i-u_{k+1}}}{\prod_{j=1, j \ne i}^{k+1} (u_j-u_i)}
+ x^{u_{k+1}-u_{k+1}} \cdot
\sum_{i=1}^{k} \frac{1}{(u_i-u_{k+1}) \prod_{j=1, j \ne i}^k (u_j-u_i)}.
\end{align*}
Thus the statement of the lemma is equivalent to the identity
\begin{equation}
\label{LagrangeEq}
\sum_{i=1}^{k} \frac{1}{(u_i-u_{k+1}) \prod_{j=1, j \ne i}^k (u_j-u_i)}=
\frac{1}{\prod_{j=1}^k (u_j-u_{k+1})}.
\end{equation}
To prove it, consider the polynomial
$$
P(x) = \sum_{i=1}^k \frac{\prod_{j=1, j \ne i}^k (u_j-x)}
{\prod_{j=1, j \ne i}^k (u_j-u_i)}.
$$
of degree $k-1$.
We have $P(u_i)=1$ for any $i \in \{ 1,2, \dots, k \}$. Hence
$P(x) \equiv 1$. The equality $P(u_{k+1})=1$ yields (\ref{LagrangeEq})
and the lemma follows.

\vspace{2mm} \noindent
{\bf Proof of Theorem \ref{MainTh3}.}
In the case $m=1$ the theorem is equivalent to Lemma 1. Now let $m>1$.
Make the substitution (\ref{Subst}) in the integral $J$ in (\ref{MainTh3Eq})
$$
J= \int_{1 \ge t_1 \ge t_2 \ge \cdots \ge t_m \ge 0}
\frac{t_m^{u_m-1}}{1-z t_m} (-\ln t_m)^s
t_1^{u_1-u_2-1} t_2^{u_2-u_3-1} \cdots t_{m-1}^{u_{m-1}-u_m-1}
dt_1 dt_2 \cdots dt_m.
$$
Applying (\ref{SimplexEq3}) for $k=m-1$ and $x=t_m$, we obtain
$$
J = \sum_{i=1}^{m} \frac{1}{\prod_{j=1, j \ne i}^{k+1} (u_j-u_i)}
\int_0^1 \frac{t_m^{u_i-1}}{1-z t_m} (-\ln t_m)^s dt_m.
$$
Use Lemma \ref{Lemma1} and the theorem follows. \\
{\bf Remark.} Theorem \ref{MainTh3} is another generalization of 
the first equality in Theorem \ref{SondowGuilleraTh1}.

\vspace{2mm}
{\bf Acknowledgment.} The author wishes to thank J.~Sondow for reading 
a preliminary version of the paper and for some useful suggestions.


\newcommand{\namefont}{\scshape}
\newcommand{\titlefont}{\itshape}

\end {document}

%% file: a4.tex
\documentclass[titlepage,12pt]{article}
\usepackage{amsfonts,amsmath,amssymb,indentfirst}

\usepackage{geometry}

%% file: library.tex
\renewcommand{\Re}{\mathop{\rm Re}\nolimits}
\def\cube#1{[0,1]^{#1}}

\newtheorem{theorem}{Theorem}
\newtheorem{lemma}{Lemma}

\newtheorem{example}{Example}
\def\proof{{\bf Proof. }}

\def\C{{\mathbb C}}

%% file: BeukersStyle.bbl
\begin{thebibliography}{99}

\bibitem{Bateman}
{\namefont Bateman H., Erdelyi A.,}
{\titlefont 
Higher Transcendental Functions
} //
Vol. 1, McGraw-Hill, New York, 1953.

\bibitem{GuilleraSondow}
{\namefont Guillera J., Sondow J.,}
{\titlefont 
Double integrals and infinite products for some classical
constants via analytic continuations of Lerch's transcendent
} //
E-print math.NT/0506319, to appear in Ramanujan J.

\bibitem{SondowMonthly}
{\namefont Sondow J.,}
{\titlefont 
Double integrals for Euler's constant and $\ln 4/\pi$ and
an analog of Hadjicostas's formula
} //
Amer. Math. Monthly 112 (2005) P. 61--65.

\end{thebibliography}
